\documentclass[12pt]{amsart}
\usepackage{amscd,amsmath,amsthm,amssymb}
\usepackage{pstcol,pst-plot,pst-3d}
\usepackage[T1]{fontenc}
\usepackage{lmodern,pst-node}
\usepackage{geometry}

\def\NZQ{\Bbb}               % the font for N,Z,Q,R,C
\def\NN{{\NZQ N}}

\def\ZZ{{\NZQ Z}}

\def\frk{\frak}               % font for "Fraktur"

\def\Phi{{\frk n}}
\def\Phi{{\frk N}}
\def\opn#1#2{\def#1{\operatorname{#2}}} % to make operators
\opn\chara{char} \opn\length{\ell} \opn\pd{pd} \opn\rk{rk}
\opn\projdim{proj\,dim} \opn\injdim{inj\,dim} \opn\rank{rank}
\opn\depth{depth} \opn\grade{grade} \opn\height{height}
\opn\embdim{emb\,dim} \opn\codim{codim}

\opn\Tr{Tr} \opn\bigrank{big\,rank}
\opn\superheight{superheight}\opn\lcm{lcm}
\opn\trdeg{tr\,deg}%\emph{
\opn\reg{reg} \opn\lreg{lreg} \opn\ini{in} \opn\lpd{lpd}
\opn\size{size} \opn\sdepth{sdepth}
\opn\link{link}\opn\fdepth{fdepth}
\opn\div{div} \opn\Div{Div} \opn\cl{cl} \opn\Cl{Cl}
\opn\Spec{Spec} \opn\Supp{Supp} \opn\supp{supp} \opn\Sing{Sing}
\opn\Ass{Ass} \opn\Min{Min}\opn\Mon{Mon}
\opn\Ann{Ann} \opn\Rad{Rad} \opn\Soc{Soc}
\opn\Im{Im} \opn\Ker{Ker} \opn\Coker{Coker} \opn\Am{Am}
\opn\Hom{Hom} \opn\Tor{Tor} \opn\Ext{Ext} \opn\End{End}
\opn\Aut{Aut} \opn\id{id}

\opn\nat{nat}
\opn\pff{pf}%   \pf exists already
\opn\Pf{Pf} \opn\GL{GL} \opn\SL{SL} \opn\mod{mod} \opn\ord{ord}
\opn\Gin{Gin} \opn\Hilb{Hilb}
\opn\aff{aff} \opn\con{conv} \opn\relint{relint} \opn\st{st}
\opn\lk{lk} \opn\cn{cn} \opn\core{core} \opn\vol{vol}
\opn\link{link} \opn\star{star}
\opn\gr{gr}

\def\pot#1#2{#1[\kern-0.28ex[#2]\kern-0.28ex]}

\opn\dirlim{\underrightarrow{\lim}}
\opn\inivlim{\underleftarrow{\lim}}
\let\union=\cup
\let\sect=\cap
\let\dirsum=\oplus

\let\Dirsum=\bigoplus

\let\to=\rightarrow

\def\Implies{\ifmmode\Longrightarrow \else
        \unskip${}\Longrightarrow{}$\ignorespaces\fi}
\def\implies{\ifmmode\Rightarrow \else
        \unskip${}\Rightarrow{}$\ignorespaces\fi}
\def\iff{\ifmmode\Longleftrightarrow \else
        \unskip${}\Longleftrightarrow{}$\ignorespaces\fi}

\let\:=\colon
\newtheorem{Theorem}{Theorem}[section]
\newtheorem{Lemma}[Theorem]{Lemma}
\newtheorem{Corollary}[Theorem]{Corollary}
\newtheorem{Proposition}[Theorem]{Proposition}

\newtheorem{Example}[Theorem]{Example}

\let\epsilon\varepsilon
\let\kappa=\varkappa
\def\qed{\ifhmode\textqed\fi
      \ifmmode\ifinner\quad\qedsymbol\else\dispqed\fi\fi}
\def\textqed{\unskip\nobreak\penalty50
       \hskip2em\hbox{}\nobreak\hfil\qedsymbol
       \parfillskip=0pt \finalhyphendemerits=0}
\def\dispqed{\rlap{\qquad\qedsymbol}}

\opn\dis{dis}
\def\pnt{{\raise0.5mm\hbox{\large\bf.}}}

\opn\Lex{Lex}

\begin{document}

\title{Stanley decompositions in localized polynomial rings}
\author{Sumiya Nasir  and Asia Rauf}
\thanks{The second author wants to acknowledge Higher Education Commission Pakistan for partial financial support during preparation of this work.}

\address{Sumiya Nasir, Abdus Salam School of Mathematical Sciences, GC University,
Lahore, Pakistan}\email{snasir@sms.edu.pk}

\address{Asia Rauf, Abdus Salam School of Mathematical Sciences, GC University,
Lahore, Pakistan} \email{asia.rauf@gmail.com}

\maketitle

\begin{abstract}
We introduce the concept of Stanley decompositions in the localized polynomial ring $S_f$ where $f$ is a product of variables, and we show that the Stanley depth does not decrease upon localization. Furthermore it is shown that for monomial ideals $J\subset I\subset S_f$ the number of Stanley spaces in a Stanley decomposition of $I/J$ is an invariant of $I/J$. For the proof of this result we introduce Hilbert series for $\ZZ^n$-graded $K$-vector spaces.
\end{abstract}

\section*{Introduction}

Let $K$ be a field, $S=K[x_1,\ldots,x_n]$ be the polynomial ring in
$n$ variables over $K$,  $A \subset \{1,2, \ldots, n\}$ and $f=\prod \limits_{j \in A}x_{j}$. Then $S_f$ is a $\ZZ^n$-graded $K$-algebra whose $K$-basis consists of the monomials  $x^{a_1}_{1}x^{a_2}_{2}\cdots x^{a_n}_{n}$ with $a_{j} \in
\ZZ$ and $a_{j}\geq 0$  if $j \notin A$.

\medskip
Standard algebraic operations like reduction module regular elements or restricted localizations have been considered in the papers \cite{R} and \cite{R2}. To better understand localization of Stanley decompositions we define in this paper Stanley decompositions of $S_f$-modules of the type $I/J$ where $J\subset I\subset S_f$ are monomial ideals in $S_f$. To do this we first have to extend the definition of Stanley spaces as it is given for $\NN^n$-graded $K$-algebras. Here the $K$-subspace of the form $uK[Z]\subset I/J$, where $u$ is a monomial in $I\setminus J$ and $Z\subset  \{x_{1}, \ldots, x_{n}\}\cup \{x_j^{-1}\: \; j\in A\}$ such that $\{x_{j}, x^{-1}_{j}\}\not\subseteq Z$ for all $j\in A$, is called {\em Stanley space} if $uK[Z]$ is a free $K[Z]$-module of $I/J$. The dimension of the Stanley space  $uK[Z]$ is defined to be $|Z|$.  As in the classical case we define a {\em Stanley decomposition} of $I/J$ as a finite direct sum of Stanley spaces ${\mathcal D}: I/J=\bigoplus^r \limits_{i=1}u_{i}K[Z_{i}],$ and set $\sdepth({\mathcal D})=\min \{|Z_{i}|: 1\leq i \leq r \}$ and $\sdepth I/J=\max \{\sdepth({\mathcal D}): {\mathcal D}\,\,  \hbox {is a Stanley decomposition of } I/J\}.$

Our more general definition of Stanley spaces is required for $S_f$, because otherwise it would be impossible to obtain a Stanley decomposition of $S_f$ when $f\neq 1$. Actually we show in Proposition~\ref{stanley} that $S_f$ has a nice canonical Stanley decomposition with Stanley spaces as defined above. Indeed, we have $S_f=  \bigoplus \limits_{L \subseteq A}
f^{-1}_{L}K[Z_{L}],$ where $$Z_{L} = \{x^{-1}_{l}| \ l \in L \} \cup
\{x_{l} | \ l \notin L \}$$ and $f_{L}=\prod \limits_{l \in L} x_{l}.$ In this decomposition all Stanley spaces are of dimension $n$ and the number of summands is $2^{|A|}$. One may ask if there exist Stanley decompositions of $S_f$ with less summands. We will show that this is not possible. Indeed, as a generalization of a result in \cite{A} we show in Theorem~\ref{spaces} that number of Stanley spaces of maximal dimension in a Stanley decomposition of an $S_f$-module of the form $I/J$ is independent of the particular Stanley decomposition of $I/J$. In order to prove this fact we introduce in Section \ref{section} a modified Hilbert function for $\ZZ^n$-graded finitely generated $K$-vector spaces $M$ with the property that $\dim_KM_a<\infty$ for all $a\in \ZZ^n$. For such modules we set $H(M, d)=\sum_{a\in \ZZ^n,\; |a|=d}\dim_K M_a$  and call the formal power series  $H_{M}(t)=\sum_{d\geq 0}H(M,d)t^d$ the Hilbert series of $M$. Here we use the notation $|a|=\sum_{i=1}^n|a_i|$ for $a=(a_1,a_2,\ldots,a_n)$. Our  Hilbert series does not as well behave as the usual Hilbert series with respect  to shifts.  Nevertheless Proposition~\ref{tricky} implies that for any monomial $u$ in $S_f$ one has  $H_{uS_f}(t)=t^{|\deg(u')|}(1+t)^{|A|}/(1-t)^n$, where $u'$ is obtained from $u$ by removing the unit factors in $u$. Quite generally we show that, just as for the usual Hilbert series,  one has  $H_{I/J}(t)=P(t)/(1-t)^d$ where $d$ is the Krull dimension of $I/J$ and $P(1)>0$.

Any monomial ideal $I\subset S_f$ is obtained by localization from a monomial ideal in $S$. Therefore, for $J\subset I\subset S$ it is natural to compare the Stanley depth of $I/J$ with that of $(I/J)_f$. Theorem~\ref{main} states that  $\sdepth I/J \leq \sdepth (I/J)_f$, and we give examples that show  that this inequality can be strict. In Theorem~\ref{fdepth} we prove a  similar result for the so-called $\fdepth$ which is define via filtrations.  Theorem~\ref{fdepth} is related to a theorem proved by the first author in \cite{S}, where it is shown that
$\sdepth S'/I_{x_1}\sect S'\geq \sdepth S/I-1$, where $S'=K[x_2,\ldots,x_n]$.

Let $J \subset I \subset S$ be monomial ideals, and consider the polynomial ring $S[t]$ over $S$. Herzog, Vladoiu and Zheng in \cite{HVZ} proved that $\sdepth(I/J)[t]=\sdepth(I/J)+1$. In Theorem~\ref{variable}, we see that a similar  relation holds  between the Stanley depth of $I/J$ and $(I/J)[t,t^{-1}]$ for monomial ideals $J\subset I\subset S_f$. It implies that Stanley depth of $I/J$, where $J\subset I\subset S_f$ monomial ideals, also if $f\neq 0$ can be computed.

\section{Monomial ideals in $S_{f}$}\label{secmon}

Let $K$ be a field and $S=K[x_1,x_2,\ldots,x_n]$ be the polynomial
ring in $n$ variables over $K.$ We fix a subset  $A \subset \{1,\ldots,n\}$  and set $$f=\prod \limits_{j \in
A}x_{j}.$$
As usual, the localization of $S$ with respect to multiplicative set $\{1, f, f^2, \ldots\}$ is denoted by $S_{f}.$
We note that
$S_{f}=K[x_{1},x_{2}, \ldots, x_{n}, \ x^{-1}_{j}:j \in A ].$

The monomials in $S_f$ are the element  $x^{a_1}_{1}x^{a_2}_{2}\cdots x^{a_n}_{n}$  with $a_{j} \in
\ZZ$ if $j \in A$ and $a_{j} \in \NN$ if $j \notin A$. The set $\Mon(S_{f})$ of monomials is a
$K$-basis of $S_{f}.$ An ideal $I \subset S_{f}$ is called a {\em monomial ideal} if it is
generated by monomials. A minimal set of monomial generators of $I$ is uniquely determined, and is denoted $G(I)$.

Any element $g \in S_{f}$ is a unique
$K$-linear combination of monomials.
\[
g=\sum a_{u}u,\quad   \text{with}\quad a_{u} \in K \quad\text{and}\quad u \in \Mon (S_{f}).
\]
We call the set $$\Supp(g)=\{ u \in \Mon (S_{f}) : a_{u}\neq 0
\}$$ the support of $g$.

For monomial $u=x^{a_1}_{1}x^{a_2}_{2}\ldots x^{a_n}_{n}$, we  set $\supp(u)=\{ i \in [n]:a_{i} \neq 0 \}$, $\supp_{+}(u)=\{i \in [n]:a_{i}>0\}$ and $\supp_{-}(u):=\{ i \in [n]:a_{i}<0\}$, and call  these sets  {\em support},  {\em positive support}, and  {\em negative support} of  $u$, respectively.

Just as for monomial ideals in the polynomial ring one shows that the set of monomials $\Mon(I)$ belonging to $I$ is a $K$-basis
of $I$. We obviously have
\begin{Proposition}
\label{monomial} Let $I \subset S_{f}$ be an ideal. The following are
equivalent.
\begin{itemize}
  \item[(a)]  $I \subset S_{f}$ is a monomial ideal.
  \item[(b)] There exists a monomial ideal $I' \subset S$  such that $I=I' S_{f}$.
\end{itemize}
If the equivalent conditions hold, then $I'$ can be chosen such that $\supp(u)\subset [n]\setminus A$ for all $u\in G(I')$.  The monomial ideal  $I'$ satisfying this extra condition is uniquely determined. Indeed, $I'=I\sect S$.
\end{Proposition}

Let $J \subset
I \subset S_{f}$ be monomial ideals. Observe that the  residue classes of the
monomials belonging to $I\setminus J$ form a $K$-basis of the residue
class $I/J.$  Therefore we may identify the classes of $I/J$ with the monomials in $I\setminus J$.

\section{Stanley decomposition of $I/J$}\label{secSta}
Let $J \subset I \subset S_{f}$ be  monomial ideals and $f=\prod \limits_{j \in
A}x_{j}$. A {\em Stanley space} of $I/J$ is a  free $K[Z]$-submodule $uK[Z]$ of $I/J$, where $u$ is a monomial of $I/J$ and
$$Z\subset \{ x_{1},\ldots,x_{n}\}\union \{x_j^{-1}\:\; j\in A\},$$
satisfying the condition that $\{x_j,x_j^{-1}\}\nsubseteq Z$ for all $j\in A$.

A {\em Stanley decomposition} of $I/J$ is  a finite direct sum of Stanley spaces
$$
{\mathcal D}: I/J=\bigoplus_{i=1}^r u_{i}K[Z_{i}].
$$

We set
$$\sdepth{\mathcal D}=\min\{|Z_{i}|: 1\leq i\leq r \},$$ and
$$\sdepth (I/J)=\max\{\sdepth {\mathcal D}:\mathcal D \text{ is  a  Stanley  decomposition of}
\ I/J\}.$$
This number is called   the {\em Stanley depth} of
$I/J.$

\begin{Proposition}
\label{stanley} The ring $S_{f}$ admits the following canonical Stanley decomposition $${\mathcal
D:S_{f}} = \bigoplus \limits_{L \subseteq A}
f^{-1}_{L}K[Z_{L}],$$ where $Z_{L} = \{x^{-1}_{l}| \ l \in L \} \cup
\{x_{l} | \ l \notin L \}$ and $f_{L}=\prod \limits_{l \in L} x_{l}.$
\end{Proposition}

\begin{proof}
Consider a monomial
$h=x^{a_{1}}_{1}\cdots x^{a_{n}}_{n}$ of $S_{f}$.  We choose  $$L=\{1\leq l \leq n:a_{l}<0 \}\subset A.$$
Then
\begin{eqnarray*}
h
  =f^{-1}_{L}(x^{a_{1}}_{1}\cdots x^{a_{n}}_{n}f_{L})
  = f^{-1}_{L}x^{b_{1}}_{1}\cdots x^{b_{n}},
\end{eqnarray*}
where $b_{l}=a_{l} \in \NN$ if $l \notin L$ and $b_{l}=a_{l}+1 \in
\ZZ$ if $l \in L$, and hence $h\in f_L^{-1}K[Z_L]$. This shows that $S_f=\sum_{L\subset A}f_LK[Z_L]$.

To show that this sum is direct, it is enough to prove any two distinct Stanley spaces in the decomposition of $I/J$ have no monomial in common, since the multigraded components of $I/J$ are of dimension $\leq 1$.

Let $L, L' \subset A$ with
$L\neq L'.$ Then  $f^{-1}_{L}\neq f^{-1}_{L'}.$ Suppose that there is a monomial $u\in I/J$ such that  $$u \in
f^{-1}_{L}K[Z_{L}]\cap f^{-1}_{L'}K[Z_{L'}],$$ that is
$$ u= f^{-1}_{L}v=f^{-1}_{L'}v'$$
for some monomials $v \in K[Z_{L}], v' \in K[Z_{L'}]$.  Since is a contradiction, because
$$L=\supp_{-}(u)=L',$$ as follows from the above equation.
\end{proof}

\begin{Example}
{\em Let $S=K[x,y,z]$   and $f=yz$. Then according to Lemma \ref{stanley},  $S_{f}=K[x,y,z]\oplus y^{-1}K[x,y^{-1},z]\oplus
z^{-1}K[x,y,z^{-1}]\oplus y^{-1}z^{-1}K[x,y^{-1},z^{-1}]$ is a Stanley decomposition of $S_{f}.$}
\end{Example}

As a consequence of Proposition~\ref{stanley} we have

\begin{Corollary}
\label{sdepth}$\sdepth (S_{f})=n.$
\end{Corollary}

\section{Localization of Stanley decompositions}
Let $J \subset I \subset S$ be monomial ideals and we set $f=\prod \limits_{j \in A}x_{j}.$ The next result shows how a Stanley decomposition of $I/J$ localizes.

\begin{Theorem} \label{main}Let ${\mathcal D}:I/J = \bigoplus^r
\limits_{i=1} u_{i}K[Z_{i}]$ be a Stanley decomposition of $I/J$.
Then ${\mathcal D_{f}}:(I/J)_{f} =\bigoplus
\limits_{i\in C}(\bigoplus \limits_{L \subset A}
 u_{i}f^{-1}_{L}K[Z^L_{i}])$ is a Stanley decomposition
of $(I/J)_{f},$ where  $C=\{i\:\; Z_A\subset Z_{i}\}$ and $Z^L_{i}=\{x^{-1}_{l} \mid l \in L,
x_{l}\in Z_{i}\}\cup \{ x_{l} \mid l \notin L, x_{l} \in Z_{i} \}$.

In particular,  we have
\[
\sdepth I/J \leq \sdepth (I/J)_{f}.
\]
\end{Theorem}

\begin{proof}
Let $u \in (I/J)_{f}$ be a nonzero monomial. We claim
that $$ u \in \sum_{i \in C} (\sum \limits_{L\subset
A}u_{i}f^{-1}_{L}K[Z^L_{i}]).
$$
Since $ u \in (I/J)_{f}$  is nonzero  it follows that $uf^a\in I\setminus J$ for all $a\gg 0$. Hence there exists an integer $i$ such that $uf^a\in u_iK[Z_i]$ for infinitely many $a>0$. This implies that $Z_A\subset Z_i$, that is, $i\in C$. So, $uf^{a} \in u_{i}K[Z_{i}]$, that is, $uf^{a}=u_{i}v$ where $v=x_1 ^{a_1}\cdots x_n ^{a_n} \in \Mon(Z_{i})$. Then $u=u_{i}vf^{-a}=u_{i}x^{b_{1}}_{1}\cdots x^{b_{n}}_{n}$, where
$$b_{i}=\left\{
  \begin{array}{ll}
     & a_{i}-a,  \,\, \hbox{if} \, \,i \in A\\
     & a_{i}, \, \, \hbox{otherwise.}
  \end{array}
\right.$$
Say $L=\supp_{-}(x^{b_{1}}_{1}\cdots x^{b_{n}}_{n})$. Note that $L \subset A$. Hence $u \in u_{i}f^{-1}_{L}K[Z^{L}_{i}]$.

\medskip
In order to prove other inclusion, consider a monomial $w \in \sum \limits_{i \in C}(\sum \limits_{L \subset A} u_{i}K[Z^{L}_{i}])$.
This implies that there exists a $i \in C$ and $L \subset A$ such that $w \in u_{i}K[Z^{L}_{i}].$ Since $Z_{A}\subset Z_{i}$, we have $wf^{a} \in u_{i}K[Z_{i}]$ for $a\gg 0.$ It follows that $wf^{a} \in I$, so $w \in I_{f}$. On the other hand,  $w \notin J_{f}$. Otherwise $wf^{a}\in J$ for $a\gg 0$, which is a contradiction, since $wf^{a} \in u_{i}K[Z_{i}]$ for $a\gg 0$.

Now we prove that the sum is direct. Let $L, H \subset A$ and
\[
u \in u_{i}f^{-1}_{L}K[Z^{L}_{i}]\cap
u_{j}f^{-1}_{H}K[Z^{H}_{j}]
\]
be a monomial such that $i\neq j$ or $L\neq H$.
Since $Z_{A} \subset Z_{i}, Z_{j}$, we have for $a\gg 0$ that
$$uf^{a} \in u_{i}K[Z_{i}]\cap
u_{j}K[Z_{j}].$$
If $i\neq j$, then $uf^{a}=0$ in $I/J$ implies $u=0$ in $(I/J)_{f}$. If $i=j$ and $L\neq H$ then $u=u_{i}f^{-1}_{L}v=u_{i}f^{-1}_{H}v'$ for some monomials $v \in K[Z^{L}_{i}]$ and $v' \in K[Z^{H}_{i}]$. Since $\supp_{-}(u_{i}f^{-1}_{L}v)=L$ and $\supp_{-}(u_{i}f^{-1}_{H}v')=H$, so $\supp_{-}(u)=L=H$, which is a contradiction.
\end{proof}

The next examples show that in the inequality on Theorem \ref{main} we may have equality or strict inequality.

\begin{Example}
{\em Let $J=(xy, yz) \subset I=(y)\subset S=K[x, y, z]$ be ideals, ${\mathcal
D}:I/J=yK[y]$ is a Stanley decomposition of $I/J.$ Thus
$\sdepth {\mathcal D}=1.$ Localizing with the
monomial $f=y,$ we obtain that ${\mathcal D_{f}}:(I/J)_{f}=yK[y]\oplus
K[y^{-1}]$ is a Stanley decomposition of $(I/J)_{f}$ and $\sdepth
{\mathcal D_{f}}=1.$}
\end{Example}
\begin{Example}
{\em Let $J=(x^3y, x^2y^2) \subset I=(x^2, xy) \subset S=K[x, y]$ be ideals, ${\mathcal
D}:I/J=xyK[y]\oplus x^2K[x]\oplus x^2yK$ is a Stanley decomposition of $I/J.$ Thus
$\sdepth {\mathcal D}=0.$ Localizing with the
monomial $f=x,$ ${\mathcal D_{f}}:(I/J)_{f}=x^2K[x]\bigoplus
xK[x^{-1}]$ is a Stanley decomposition of $(I/J)_{f}$ and $\sdepth
{\mathcal D_{f}}=1.$}
\end{Example}

The following example shows that a Stanley decomposition of $I/J$ which gives the Stanley depth of $I/J$  may localize to a Stanley decomposition whose Stanley depth is smaller than the Stanley depth of the localization of $I/J$.
\begin{Example}\label{mom}
{\em Let $I=(x, y, z) \subset S=K[x, y, z]$ be the graded maximal ideal of $S$

${\mathcal
D}:I=xK[x, y]\oplus yK[y, z]\oplus zK[x, z]\oplus xyzK[x, y, z]$ is a Stanley decomposition of $I$. Thus
$\sdepth {\mathcal D}=2$ which is also the Stanley depth of $I$.  Localizing with the
monomial $f=x,$ we get the Stanley decomposition ${\mathcal D}_f$ of $I_f$ which is
\[
xK[x, y]\oplus K[x^{-1}, y]\oplus zK[x, z]
\oplus zx^{-1}K[x^{-1}, z]\oplus xyzK[x, y, z]\oplus yzK[x^{-1}, y, z].
\]
Thus $\sdepth{\mathcal D_{f}}=2.$ However  $I_{f}=K[x^{\pm 1}, y, z]$, and  hence  $\sdepth I_{f}=3.$}
\end{Example}

\section{Stanley decompositions and polynomial extensions}

Herzog et al. in \cite[Lemma 3.6]{HVZ} proved that the Stanley depth of the module  increases by one in a polynomial ring extension by one variable. We generalize this result to  localized rings.
\begin{Theorem}
\label{variable}Let $J \subset I \subset S_{f}$ be monomial ideals. Then
\[
\sdepth (I/J)[t]=\sdepth (I/J)[t,t^{-1}]=\sdepth I/J+1.
\]
\end{Theorem}
\begin{proof}
Let ${\mathcal D}: I/J  = \bigoplus^r \limits_{i=1}
v_{i}K[Z_{i}]$ be a Stanley decomposition of $I/J$.  Then we obtain the Stanley decompositions
\[
(I/J)[t]=\bigoplus^r \limits_{i=1} v_{i}K[Z_{i}][t]
=\bigoplus^r \limits_{i=1} v_{i}K[Z_{i},t]
\]

\[
 ((I/J)[t])_{t}=\bigoplus^r \limits_{i=1} (v_{i}K[Z_{i},t]\oplus
v_{i}t^{-1}K[Z_{i},t^{-1}]).\]
 This implies that
\[1+ \sdepth I/J \leq
\sdepth (I/J)[t,t^{-1}],\;  \sdepth (I/J)[t].
\]
Let ${\mathcal D'}:(I/J)[t,t^{-1}]=\bigoplus^r \limits_{i=1}
v_{i}K[W_{i}]$ be a Stanley decomposition of
$(I/J)[t,t^{{-1}}]$. Next we show  that
\[
I/J=\bigoplus_{i \in [r]} v_{i}K[W_{i}]\cap S_{f},
\]
and that each $v_{i}K[W_{i}]\cap S_{f}=u_{i}K[W_{i}\backslash \{ t,t^{-1}\}]$ for a suitable monomial $u_i\in S_{f}$, provided  $v_{i}K[W_{i}]\cap S_{f}\neq 0$.

Since the direct sum $\mathcal{D}'$  commutes with taking the intersection with $S_f$ and since $(I/J)[t,t^{{-1}}]\sect  S_{f}=I/J$, the desired decomposition of $I/J$ follows.

Suppose $v_{i}K[W_{i}]\cap S_{f}\neq 0$. Then there exists  monomials $v\in S_{f}$ and $w\in K[W_i]$ such that $v=v_iw$. We may assume that $v_i$ does not contain $t$ as a factor, because $t^{-1}$ must be a factor of $w$ which implies that $t^{-1}\in W_i$. Thus we may replace $v_i$ by the monomial which is obtained from $v_i$ by removing the power of $t$ which appears in $v_i$. Similarly we may assume $t^{-1}$ is not a factor of $v_i$. Then it follows that  $v_{i}K[W_{i}]\cap S_{f}$ consists of all monomials $v_iw$ with $w\in K[W_i]$ where neither $t$ nor $t^{-1}$ is a factor of $w$. In other words, $w\in K[W_i\setminus\{t,t^{-1}\}]$.

From these considerations we conclude that $\ 1 + \sdepth I/J \geq
\sdepth (I/J)[t,t^{-1}]$. In the same way one proves the inequality  $\ 1 + \sdepth I/J \geq
\sdepth (I/J)[t]$. This yields the desired conclusions.
\end{proof}

An immediate consequence of Theorem  \ref{variable} is the following
\begin{Corollary}
\label{rabia}
\[
\sdepth (I/J)[t^{\pm 1}_{1}, \ldots, t^{\pm 1}_{r}]=\sdepth I/J+r.
\]
\end{Corollary}

Let $J\subset I\subset S_f$ be monomial ideals, and  $S'=K[x_i:i\not\in A]\subset S$.  Then there exist monomial ideals $J'\subset I'\subset S'$  such that $J'S_f=J$ and $I'S_f=I.$ We have
\[
I'/J'[x_i ^{\pm1}:i\in A]=I/J.
\]
Hence $\sdepth I'/J'+|A|=\sdepth I/J$, by Corollary \ref{rabia}. Since the Stanley depth of $I'/J'$ can be computed in finite number of steps by the method given by Herzog, Vladoiu and Zheng \cite{HVZ}, the  Stanley depth of $I/J$ can  be computed as well in a finite number of steps.

\section{Fdepth and localization}

Let $J\subset I\subset S$ be monomial ideals of $S$.
A chain  ${\mathcal F}:J=J_{0}\subset J_{1}\subset \ldots \subset J_{r}=I$ of monomial ideals is called  a prime filtration of $I/J$, if $J_{i}/J_{i-1}\cong S/P_{i}(-a_{i})$ where $a_{i} \in \ZZ^n$ and where each $P_{i}$ is a monomial prime ideal. We call the set of prime ideals $\{P_{1}, \ldots, P_{r}\}$ the support of ${\mathcal F}$ and denote it by $\Supp({\mathcal F}).$ In \cite{HVZ}, Herzog, Vladoiu and Zheng define $\fdepth$ of $I/J$ as follows:
 \[
 \fdepth ({\mathcal F})=\min \{ \dim S/P : P \in \Supp({\mathcal F}) \}
  \] and
\[
\fdepth I/J= \max \{ \fdepth {\mathcal F} : {\mathcal F} \, \text{is a prime filtration of} \, I/J \}.
\]
These definitions can be immediately extended to monomial ideals $J\subset I\subset S_{f}$. We then get

\begin{Theorem} \label{fdepth}
\[
\fdepth I/J \leq \fdepth (I/J)_{f}.
\]
\end{Theorem}
\begin{proof}
Let ${\mathcal F}\: J=J_{0}\subset J_{1}\subset\ldots \subset J_{r}=I$  be a prime filtration of $I/J$ with  $\Supp({\mathcal F})=\{P_1,\ldots,P_r\}$. Then we obtain a filtration ${\mathcal F}_{f}\: J_{f}=(J_{0})_{f}\subset (J_{1})_{f}\subset \ldots \subset (J_{r})_{f}=I_{f}$ with factors $S_{f}/P_iS_{f}$. Thus if we drop the ideals $(J_i)_{f}$ for which $S_{f}/P_iS_{f}=0$ we obtain a prime filtration of $(I/J)_{f}$. Since $\dim S_{f}/P_iS_{f}=\dim S/P_i$, whenever $\dim S_{f}/P_iS_{f}\neq 0$, it follows that $\fdepth {\mathcal F}_{f}\geq \fdepth {\mathcal F}$. This implies the desired conclusion.
\end{proof}

\section{Hilbert series of multigraded $K$-vector spaces}
\label{section}
Let $J\subset I\subset S_f$ be monomial ideals. In this section we want to show that the number of maximal Stanley spaces in any Stanley
decomposition of $I/J$ is an invariant of $I/J$. To prove this we introduce a new kind of Hilbert series.

Note that the localized ring $S_{f}$ is naturally $\ZZ^{n}$-graded  with $\ZZ^{n}$-graded components
$$(S_f)_{a}=\left\{
  \begin{array}{ll}
    Kx^{a}, &  \hbox{} a_{i}\geq 0\, \, \hbox{if} \, \,i \notin \supp(f);\\
     0,&   \hbox{otherwise.}
  \end{array}
\right.$$

Let $M=\bigoplus \limits_{a \in \ZZ^{n}}M_{a}$  be $\ZZ^{n}$- graded $K$-vector space with $\dim_K M_a<\infty$ for all $a\in\ZZ^n$.  Then for all $d \in \NN$ we define
\[
M_{d}=\bigoplus \limits_{|a|=d}M_{a}\quad  \text{where}\quad |a|=\sum^r \limits_{j=1}|a_{j}|,
\]
and set  $H(M, d)=\dim_{k}M_{d}.$ We call the function  $H(M,-)\: \NN \to \NN$ the {\em Hilbert function}, and the  series $H_{M}(t)=\sum \limits_{d\geq 0}H(M, d)t^{d}$ the {\em Hilbert series} of $M$.

We consider an example to illustrate our definition. Figure 1 displays the ideal $I=(x^3,x^2y,y^2) \subset S=K[x, y]$.  The grey area represents the monomial $K$-vector space spanned by the monomials in $I$. The slanted lines represent $S_{d}={\bigoplus_ {a \in \NN^{2}\atop |a|=d}}S_{a}$

\bigskip
%%%%%%%%%%%%%%%%%%%%%%%%%%%%
\begin{figure}[hbt]
\begin{center}
\begin{pspicture}(-1,-1)(6,6)
%\psgrid[subgriddiv=1,griddots=10,gridlabels=10pt](0,0)(6,6)
%\psdot*(1,5)
\pspolygon*[linecolor=lightgray,fillstyle=solid, fillcolor=lightgray](0,2)(0,5)(5,5)(5,0)(3,0)(3,1)(2,1)(2,2)(0,2)
%\pspolygon*[linecolor=lightgray,fillstyle=solid, fillcolor=lightgray](0,2)(0,5)(2,5)(2,2)
%\pspolygon*[linecolor=lightgray,fillstyle=solid, fillcolor=lightgray](2,1)(2,5)(3,5)(3,1)
\psline[arrowsize=6pt]{->}(0,0)(0,6)\psline[arrowsize=6pt]{->}(0,0)(6,0)
\psline[linestyle=dashed](1,0)(0,1)
\psline[linestyle=dashed](2,0)(0,2)
\psline[linestyle=dashed](3,0)(1,2)\psline(1,2)(0,3)
\psline(0,4)(4,0)
\psline(0,5)(5,0)
%%%%%%%%%%%%%%%%%%%%%%Dots%%%%%%%%%%%%%%%%%%%
\psdot*(0,0)\psdot*(0,1)\psdot*(0,2)\psdot*(0,3)\psdot*(0,4)\psdot*(0,5)
\psdot*(1,0)\psdot*(1,1)\psdot*(1,2)\psdot*(1,3)\psdot*(1,4)\psdot*(1,5)
\psdot*(2,0)\psdot*(2,1)\psdot*(2,2)\psdot*(2,3)\psdot*(2,4)\psdot*(2,5)
\psdot*(3,0)\psdot*(3,1)\psdot*(3,2)\psdot*(3,3)\psdot*(3,4)\psdot*(3,5)
\psdot*(4,0)\psdot*(4,1)\psdot*(4,2)\psdot*(4,3)\psdot*(4,4)\psdot*(4,5)
\psdot*(5,0)\psdot*(5,1)\psdot*(5,2)\psdot*(5,3)\psdot*(5,4)\psdot*(5,5)
%%%%%%%%%%%%%%%%%%%%% End Dots %%%%%%%%%%%%%%%%%%%%%%%%
\psline[linewidth=1.2pt](3,0)(3,1)\psline[linewidth=1.2pt](3,1)(2,1)\psline[linewidth=1.2pt](2,1)(2,2)\psline[linewidth=1.2pt](2,2)(0,2)
%%%%%%%%%%%%%%%%%%%%% Text %%%%%%%%%%%%%%%%%%%%%%%
\rput(1,-0.4){$x$}\rput(2,-0.4){$x^2$}\rput(3,-0.4){$x^3$}
\rput(-0.4,1){$y$}\rput(-0.4,2){$y^2$}\rput(-0.4,3){$y^3$}
\end{pspicture}
\end{center}
\caption{}\label{Fig1}
\end{figure}

%%%%%%%%%%%%%%%%%%%%%%%%%%%%%%%%%%%%%%%%%%%%%%%%%%%%%%%%%%%%%%
The $S$-modules $I(a)$ with $a=(4,3)$ is an $S$-submodule  of $S_{f}$ where $f=xy$. In Figure 2, the grey area displays $I(a)$. For all $d\in \NN$, the boundaries of the squares of diagonal length $2d$, represent $\bigoplus_{|a|=d}(S_{f})_{a}$.
%%%%%%%%%%%%%%%%%%%%%%%%%%%%%%%%%%%%%%%%%%%%%%%%%%%%%%%%%%%%%%%%%%%%
\begin{figure}[hbt]
\begin{center}
\begin{pspicture}(-4,-4)(4,4)
\pspolygon*[linecolor=lightgray,fillstyle=solid, fillcolor=lightgray](-2,0)(-2, 3.5)(3.5,3.5)(3.5,-1.5)(1.5,-1.5)(-0.5,-1.5)(-0.5,-1)(-1,-1)(-1,-0.5)(-2,-0.5)(-2,0)
%(-2,1)(0,3)(3,0)(1.5,-1.5)(-0.5,-1.5)(-0.5,-1)(-1,-1)(-1,-0.5)(-2,-0.5)(-2,0)
\psline[arrowsize=8pt, linewidth=1.5pt]{<->}(-4,0)(4,0)
\psline[arrowsize=8pt, linewidth=1.5pt]{<->}(0,-4)(0,4)
\pspolygon(-0.5,0)(0,0.5)(0.5,0)(0,-0.5)
\pspolygon(-1,0)(0,1)(1,0)(0,-1)
\pspolygon(-1.5,0)(0,1.5)(1.5,0)(0,-1.5)
%%%%%%%%%%%%%%%%%%%%%%% Polygon & Dots %%%%%%%%%%%%%%%%%%%%%%%%%%%%%%%%%%%%%%%%%%%%
\psline(-2,0)(0,2)\psline(0,2)(2,0)\psline(2,0)(0,-2)\psline[linestyle=dashed](0,-2)(-1.5,-0.5)\psline(-1.5,-0.5)(-2,0)
\psline[linestyle=dashed](-2.5,0)(-2,0.5)\psline(-2,0.5)(0,2.5)\psline(0,2.5)(2.5,0)\psline(2.5,0)(1,-1.5)\psline[linestyle=dashed](1,-1.5)(0,-2.5)\psline[linestyle=dashed](0,-2.5)(-2.5,0)
\psline[linestyle=dashed](-3,0)(-2,1)\psline(-2,1)(0,3)\psline(0,3)(3,0)\psline(3,0)(1.5,-1.5)\psline[linestyle=dashed](1.5,-1.5)(0,-3)\psline[linestyle=dashed](0,-3)(-3,0)
\psline(4,-1.5)(-0.5,-1.5)\psline(-0.5,-1.5)(-0.5,-1)\psline(-0.5,-1)(-1,-1)\psline(-1,-1)(-1,-0.5)\psline(-1,-0.5)(-2,-0.5)\psline(-2,-0.5)(-2,4)
\psdot*[dotsize=4pt](0,-1.5)\psdot*(-0.5,-1.5)\psdot*(-0.5,-1)\psdot*(-1,-1)\psdot*(-1,-0.5)\psdot*(-1.5,-0.5)\psdot*(-2,-0.5)\psdot*[dotsize=4pt](-2,0)
\end{pspicture}
\end{center}
\caption{}\label{Fig2}
\end{figure}
%%%%%%%%%%%%%%%%%%%%%%%%%%%%%%%%%%%%%%%%%%%%%%%%%%%%%%%%%%%%%%%%%%%%%%%%%%%%

In this example, if we consider $I(a)$ as a graded $S$-module in the usual sense, then $I(a)_4=I_{11}$, and hence $\dim_KI(a)_4=12$. If we apply our new definition of $M_d$ to $I(a)$, then  $\dim_KI(a)_4=15$, as can be seen in the picture.

\medskip
In case that $f=1$, and $M$ is a finitely generated $\ZZ^n$-graded $S$-module with the property that the multidegree of all generators of $M$ belong to $\NN^n$, then  our definition of the Hilbert function coincides with the usual definition. For properties of  classical Hilbert series we refer to  the book \cite{BH}.

If we would define  $H(M,d)$ in the usual way as $\dim_K(\Dirsum_{a,\; \sum_ia_i=d}M_a)$, then this number would be in general infinite. On the other hand, our definition has the drawback, that the components $(S_f)_d$ don't define a grading of $S_f$. In other words,  we do not have in general that $(S_f)_{d_{1}}(S_f)_{d_{2}}\subset (S_f)_{d_{1}+d_{2}}$ for all $d_{1}$ and $d_{2}$. Nevertheless, we shall prove that $H(I/J, d)$ is a polynomial function for
$d\gg 0$.

In the case that $J \subset I \subset S_{f}$ are  monomial ideals,  then the Hilbert function of $I/J$ is given by $$H(I/J, d)= |\{ a \in \ZZ^{n}| \, x^{a} \in I \setminus J \, \hbox{and} \, |a|=d \}|.$$

Let $u=x^{a_{1}}_{1} \cdots x^{a_{n}}_{n} \in \Mon(S_{[n]})$. We set  $|u|=x^{|a_{1}|}_{1} \cdots x^{|a_{n}|}_{n}$. Then $|u|\in \Mon(S)$.

Let $S_{[n]}=K[x^{\pm1}_{1},\ldots, x^{\pm1}_{n}]$ and $K[Z]\subset S_{[n]}$ be a Stanley space. We set $$\overline{Z}= \{x_{i}:x_{i}\in Z \text{ or } x^{-1}_{i}\in Z\}.$$
\begin{Lemma}
\label{iso}
Let $v\in \Mon(S_{[n]})$ such that
\[
\supp_{+}(v)\cap\{i:{x^{-1}_{i}}\in Z\}=\emptyset
\] and
\[
\supp_{-}(v)\cap\{i:{x_{i}}\in Z\}=\emptyset.
\]
Then $vK[Z]\cong |v|K[\overline{Z}] $, and for each $d \in \NN$ we have $$(vK[Z])_{d}\simeq (|v|K[\overline{Z}])_{d}.$$
\begin{proof}

Since $v\in \Mon(S_{[n]})$ such that $\supp_{+}(v)\cap\{i:{x^{-1}_{i}}\in Z\}=\emptyset$ and $\supp_{-}(v)\cap\{i:{x_{i}}\in Z\}=\emptyset$,
it follows that $|vw|=|v||w|$ for all $w\in K[Z]$. Therefore the map $u\mapsto |u|$ induces for each $d$ an isomorphism $(vK[Z])_{d}\simeq  (|v|K[\overline{Z}])_{d}$ of $K$-vector spaces.
\end{proof}
\end{Lemma}

\begin{Proposition}
\label{tricky}
Let $T_{A}=K[x^{\pm1}_{i}:i \in A, x_{i_{1}},\ldots, x_{i_{t}}]\subset S_{[n]}$, where $\{i_{1},\ldots,i_{t}\}\subset [n]\backslash A$.
Let $u\in\Mon(S_{[n]})$, $u=x_1^{a_1}x_2^{a_2}\cdots x_n^{a_n}$ such that $a_{i}\geq 0$ for $i \notin A$. We set $u'=\prod_{i\not\in A}x_i^{a_i}$. Then
\[
H_{uT_{A}}(t)=H_{u'T_{A}}(t)=t^{|\deg(u')|}H_{T_{A}}(t)=t^{|\deg(u')|}\frac{(1+t)^{|A|}}{(1-t)^{d}},
\]
where $d=\dim T_{A}$.
\end{Proposition}
\begin{proof}  Since $uT_A=u'T_A$, it follows that $H_{uT_A}(t)=H_{u'T_A}(t)$. In order to prove the second equality,  we use the Stanley decomposition $T_A=\Dirsum_{L\subset A}f_L^{-1}K[Z_L]$ where $f_{L}=\prod_{i \in L}x_{i}$ and $Z_{L}=\{x^{-1}_{i}:i \in L\}\cup\{x_{i}:i \in A\backslash L\}\cup\{x_{i_{1}}, \ldots, x_{i_{t}}\}$, from which we obtain the Stanley decomposition
\[
u'T_A=\Dirsum_{L\subset A}u'f_L^{-1}K[Z_L].
\]
Applying Lemma~\ref{iso}, we obtain from this decomposition the following identities
\begin{eqnarray*}
H_{u'T_A}(t) &=&\sum_{k=0}^{|A|} \sum_{L\subset A\atop |L|=k} H_{u'f_L^{-1}K[Z_L]}(t)= \sum_{k=0}^{|A|} \sum_{L\subset A\atop |L|=k} H_{u'f_LK[\overline{Z_{L}}]}(t)\\
&=& \sum_{k=0}^{|A|} \sum_{L\subset A\atop |L|=k}\frac{t^{k+|\deg(u')|}}{(1-t)^d}=t^{|\deg(u')|}\sum_{k=0}^{|A|}{|A|\choose k}\frac{t^{k}}{(1-t)^d}\\
&=& t^{|\deg(u')|}\frac{(1+t)^{|A|}}{(1-t)^d}=t^{|\deg(u')|}H_{T_A}(t),
\end{eqnarray*}
since $T_A=\Dirsum_{L\subset A}f_L^{-1}K[Z_L]$ and $H_{f^{-1}_{L}K[Z_{L}]}= \frac{t^{|L|}}{(1-t)^d}$ by Lemma \ref{iso}, we get that
\[
H_{T_{A}}(t)=\sum_{L \subset A}H_{f^{-1}_{L}K[Z_{L}]}(t)=\frac{{(1+t)}^{|A|}}{(1-t)^{d}}.
\]
\end{proof}
It follows from the next result that the Hilbert series of finitely generated $\ZZ^{n}$-graded $S_{f}$-module $I/J$ can be written as a rational function with denominator $(1-t)^d$, where $d=\dim I/J$.
\begin{Theorem}
Let $A\subset [n]$,  $f=\prod_{i\in A}x_i$, and $J\subset I\subset S_f$ be monomial ideals. Then $I/J$ admits a filtration
\[
0=M_0\subset M_1\subset\ldots\subset M_r=I/J
\]
of $\ZZ^n$-graded $S_f$-submodules of $I/J$ such that for each $i$ we have $$M_{i+1}/M_i\cong(S_f/P_i)(-a_i),$$ where $\{P_1,\ldots,P_r\}$ is a set of $\ZZ^n$-graded prime ideals of $S_f$ containing all minimal prime ideals of $I/J$, and where $a_i\in\NN^n$ with $a_i(j)=0$ if $j\in A.$ Moreover,
\[
H_{I/J}(t)=\sum_{i=1}^rH_{S_f/P_i(-a_i)}(t)\text{ and } H_{S_f/P_i(-a_i)}(t)=\frac{Q_i(t)(1+t)^{|A|}}{(1-t)^{d_i}},
\]
where $d_i=\dim S_f/P_i$ and $Q_i(t)$ is a polynomial with $Q_i(1)=1.$ In particular, $H_{I/J}(t)=Q(t)/(1-t)^d$ with $Q(1)=k\cdot 2^{|A|}$, where
\[
k=|\{i: \dim S_f/P_i=d\}|\quad \text{and}\quad d=\dim I/J.
\]
\end{Theorem}
\begin{proof}
Since $J\subset I\subset S_f$ monomial ideals, we may assume that $\supp(u)\subset [n]\setminus A$ for all $u\in G(I)\cup G(J)$, see  Proposition \ref{monomial}. Let $S'=K[x_i:i\not\in A]\subset S$ be the polynomial ring over $K$. Then there exist monomial ideals $J'\subset I'\subset S'$ such that $J=J'S_f$, $I=I'S_f$. Consider a prime filtration $$J'=I_{0} \subset I_{1}\ldots \subset I_{r}=I'$$ of $I'/J'$ where $I_{i+1}/I_{i}\cong (S'/p_{i})(-a_{i}),$ where $p_{i}$ is a monomial prime ideal  $a_i\in\NN^n$ with $a_i(j)=0$ if $j\in A.$  It follows that $$J=J'S_{f}=I_{0}S_{f} \subset I_{1}S_{f}\ldots \subset I_{r}S_{f}=I$$ is a prime filtration of $I/J$ with  $I_{i+1}S_f/I_iS_f\cong S_f/P_i(-a_i)$ where $p_iS_f=P_i$.
Because of the additivity of our Hilbert function, it suffices to show that
\[
H_{S_f/P_i(-a_i)}(t)=Q_i(t)H_{S_f/P_i}(t)
\]
for some polynomial $Q_i(t)$ with $Q_i(1)\neq0$. But this follows from Proposition \ref{tricky}. The rest of the statements are obvious.
\end{proof}

For the proof of our main result, we need the following
\begin{Lemma}
\label{importamt}
Let $uK[Z]$ be a Stanley space of $S_f$. Then $H_{uK[Z]}(t)=Q(t)/(1-t)^m$ where  $Q(1)=1$ and $m=|Z|$.
\end{Lemma}

\begin{proof}
Let $u=x_1^{a_1}\cdots x_n^{a_n}$, and let $\mathcal{I}$ be the set of indices $i$ for which either $a_i>0$ and $x_i^{-1}\in Z$,  or $a_i<0$ and $x_i\in Z$. Then we set $r=\sum_{i\in \mathcal{I}}|a_i|$, and prove  the lemma by induction on  $s=\min\{r,m\}$. If $m=0$, then the assertion is trivial, and if $r=0$, then the result  follows from Lemma~\ref{iso}.

Now let $s>0$. Then we may assume that $a_1>0$ and $x_1^{-1}\in Z$. In this case we  have the following direct sum of $K$-vector spaces
\[
uK[Z]=uK[Z\setminus\{x_1^{-1}\}]\dirsum vK[Z],
\]
where $v=x_1^{a_1-1}\cdots x_n^{a_n}$. In the first summand $m$ is smaller and in the second $r$ is smaller than the corresponding numbers for $uK[Z]$. Thus we may apply our induction hypothesis according to which there exist polynomials $Q_1(t)$ and $Q_2(t)$ with $Q_i(1)=1$, and  such that
\[
H_{uK[Z\setminus\{x_1^{-1}\}]}(t)= Q_1(t)/(1-t)^{m-1}\quad \text{and}\quad H_{vK[Z]}(t)= Q_2(t)/(1-t)^{m}.
\]
It follows that $H_{uK[Z]}(t)=Q(t)/(1-t)^m$ with $Q(t)=Q_1(t)(1-t)+Q_2(t)$. Since $Q(1)=1$, the proof is completed.
\end{proof}

\begin{Theorem}
\label{spaces}Let  $J\subset I\subset S_f$ be monomial ideals,   $\mathcal{D}: I/J=\bigoplus^r \limits_{i=1}u_{i}K[Z_{i}]$  a Stanley decomposition of $I/J$ and $d=\max\{|Z_{i}|: i\in 1,2,\ldots,r\}$. Then $H_{I/J}(t)=P_{I/J}(t)/(1-t)^d$ with $P_{I/J}(1)= |\{i\:\; |Z_i|=d\}|$.

\end{Theorem}
\begin{proof}
We have
\[
H_{I/J}(t)=\sum^r_{i=1}H_{u_{i}K[Z_{i}]}(t)
\]

By Lemma \ref{importamt}, for all $i\in\{1, 2,\ldots, r\}$ we obtain that $H_{u_{i}K[Z_{i}]}(t)=\frac{Q_{i}(t)}{(1-t)^{|Z_{i}|}}$, where $Q_{i}(1)=1$. Thus
\[
H_{I/J}(t)=\sum^r \limits_{i=1}\frac{Q_{i}(t)}{{(1-t)^{|Z_{i}|}}}=\frac{P_{I/J}(t)}{(1-t)^{d}}
\] where,
$P_{I/J}(t)=\sum^r \limits_{i=1}(1-t)^{d-|Z_{i}|}Q_{i}(t)$.
It follows that $P_{I/J}(1)= |\{i\:\; |Z_i|=d\}|$.
\end{proof}

This Theorem implies that the number of Stanley spaces of maximal dimension in a Stanley decomposition of $I/J$ is independent of the special Stanley decomposition of $I/J$.

Proposition~\ref{stanley} and Theorem~\ref{spaces}  yield the following result.

\begin{Corollary}
\label{maximal}
The number of Stanley spaces of maximal
dimension in any Stanley decomposition of $S_{f}$ is equal to
$2^{|A|}.$
\end{Corollary}
We conclude our paper with the following concrete example.

\begin{Example}
{\em
Let $S=K[x, y, z]$ and $f=z$. Then $S_{f}=K[x, y, z^{\pm1}]$.
Let $J=(x^{2})\subset I=(x, y^2)\subset S_{f}$ be monomial ideals in $S_f$. A Stanley decomposition of $I/J$ is $xK[y, z]\oplus xz^{-1}K[y]\oplus xz^{-2}K[y, z^{-1}]\oplus y^2K[y, z]\oplus y^2z^{-1}K[y, z^{-1}]$. Thus in any other Stanley decomposition of $I/J$ the number of maximal Stanley spaces is 4. Calculating the Hilbert function of $I/J$ by using this Stanley decomposition we find that  $H_{I/J}(t)=\frac{t+2t^2+t^3}{(1-t)^{2}}$. Thus $P(t)= t+2t^2+t^3$, and $P(1)=4$, as expected. }
\end{Example}

\end{document}